\documentclass[11pt]{article}
\usepackage{amsmath}
\usepackage{amsthm}
\usepackage{fullpage}
\usepackage{graphicx}
\usepackage{sublabel}
\newcommand {\norm} [1] { \lVert #1 \rVert}

\newcommand {\abs} [1] {\left| #1 \right|}
\newcommand {\Abs} [1] {\bigl\lvert #1 \bigr\rvert}

\newcommand {\E}{{\mathrm{E}}}
\renewcommand{\P}{{\mathrm{P}}}
\newcommand{\var}{\mathop{\mathrm{Var}}}
\newcommand{\covar}{\mathop{\mathrm{Covar}}}
\newcommand{\des}{\mathop{\mathrm{des}}}
\newcommand{\inv}{\mathop{\mathrm{inv}}}

\newtheorem {thm} {Theorem}[section]

\newtheorem {lem} [thm] {Lemma}

\theoremstyle {definition}

\theoremstyle {remark}

\numberwithin{equation}{section}
\title{Normal Approximations for Descents and Inversions of Permutations of Multisets}
\author{Mark Conger and D. Viswanath \thanks{ Department of Mathematics, University of Michigan, 
530 Church Street, Ann Arbor, MI 48109, U.S.A.,
{\tt mconger@umich.edu and divakar@umich.edu}.
This work was supported by a research fellowship from the Sloan Foundation.}}
\begin{document}
\maketitle
\begin{abstract}
Normal approximations for descents and inversions of permutations of the set
$\{1,2,\ldots,n\}$ are well known.
We consider the number of
inversions of a permutation $\pi(1), \pi(2),\ldots,\pi(n)$ of a multiset with 
$n$ elements, which is the number of pairs $(i,j)$ with $1\leq i < j \leq n$
and $\pi(i)>\pi(j)$. The number of descents is the number of $i$ in the range
$1\leq i < n$ such that $\pi(i) > \pi(i+1)$. We prove that, appropriately normalized,
the distribution of both inversions and descents of a random permutation of the multiset
approaches the normal distribution as $n\rightarrow\infty$, provided that the permutation
is equally likely to be any possible permutation of the multiset and  no element occurs
more than $\alpha n$ times in the multiset for a fixed $\alpha$ with
$0<\alpha < 1$. 
Both normal approximation theorems are proved using the size bias version of Stein's
method of auxiliary randomization and are accompanied by error bounds.
\end{abstract}

\section{Introduction}
Let $\pi(1), \pi(2), \ldots, \pi(n)$ be a permutation of the multiset
$\{1^{n_1},2^{n_2},\ldots, h^{n_h}\}$ with
$n_1+\cdots+n_h = n$.  The number of inversions,
denoted $\inv(\pi)$, is defined as the number of pairs $(i,j)$ with
$1\leq i < j \leq n$ and $\pi(i) > \pi(j)$.  The number of descents,
denoted $\des(\pi)$, is the number of positions $i$ with $1\leq i < n$
and $\pi(i) > \pi(i+1)$. Assume that $\pi$ is uniformly distributed. In
this article, we use Stein's method to prove normal approximations with
error bounds for $\inv(\pi)$ and $\des(\pi)$.

In the special case where $\pi$ is a uniformly distributed permutation
of the set $\{1,2,\ldots,n\}$, the distributions of both $\inv(\pi)$ and
$\des(\pi)$ admit simple descriptions. The distribution of $\inv(\pi)$
is equal to that of the sum $X_1+\cdots+X_{n-1}$, where the random variables
$X_i$, $1\leq i\leq n-1$, are  independent with $X_i$ uniformly distributed
over the set $\{0,1,\ldots, i\}$. To obtain the distribution of $\des(\pi)$,
we need the sum $X_1+\cdots+X_{n-1}+X_n$, where the $X_i$ are
independent and uniformly distributed in the interval $[0,1]$. 
The probability that this sum lies in the interval $[d,d+1]$ equals
the probability that $\des(\pi)$ equals $d$. According to Knuth 
\cite{KnuthVol3}, the first of these two results was noticed by
O. Rodriguez in 1839.  The result about $\des(\pi)$ was alluded to
by Barton and Mallows \cite{BM65}. An elegant proof is due to Stanley
\cite{Stanley1}.

Normal approximations to $\des(\pi)$ and $\inv(\pi)$ in this special 
case can be obtained using these results and standard versions of the
central limit theorem.
The bounds
\begin{align}
\sublabon{equation}
\Bigl|P\Biggl(\frac{\des(\pi) - (n-1)/2}{\sqrt{(n+1)/12}}\leq x\Biggr)-\Phi(x)\Bigr|
&\leq \frac{C}{\sqrt{n}} \label{eqn-1-1a}\\
\Bigl|P\Biggl(\frac{\inv(\pi)-\frac{1}{2}\binom{n}{2}}{\sqrt{n(n-1)(2n+5)/72}}\leq x
\Biggr)-\Phi(x)\Bigr| &\leq \frac{C}{\sqrt{n}},
\label{eqn-1-1b}
\end{align}
\sublaboff{equation}
\noindent where $C$ is a constant and $\Phi$ is the standard normal distribution, were proved
using the method of exchangeable pairs \cite{RR97, SteinBook} 
by Fulman \cite{Fulman04}. Other proofs of \eqref{eqn-1-1a}
using Stein's method are sketched
in \cite{Conger1} and \cite{Goldstein04}.

 From the survey by Barton and Mallows \cite{BM65}, it appears
that the asymptotic normality of a quantity closely related to
$\des(\pi)$, where $\pi$ is a uniformly distributed permutation of
the set $\{1,\ldots,n\}$, was stated by Bienyam\'e in 1874
(Bull. Soc. Math. France, vol. 2, p. 153-154). Bienyam\'e was interested
in statistical applications. So were Levene and Wolfowitz \cite{LW1} who stated
that {\it runs} were widely used in quality control and in the study of
economic time series. Runs are the monotone segments within a sequence
of numbers and are closely related to descents. 
An early proof of the asymptotic normality of descents, which is
implied by \eqref{eqn-1-1a}, is due to Wolfowitz \cite{Wolfowitz44}.

Let $\{1^{n_1}, 2^{n_2},\ldots, h^{n_h}\}$ be a multiset, where $n_a$, $1\leq a\leq h$, are positive integers. Let $n=n_1+n_2+\ldots+n_h$ be the number of elements of the multiset. 
Let $\alpha$ be a fixed number in $(0,1)$.
We assume that $n_a \leq \alpha n$ for $1\leq a\leq h$.
Let $\pi$ be a uniformly distributed permutation of this multiset.
We consider $\inv(\pi)$ and $\des(\pi)$ in this more general situation.
The bounds
that we obtain for the errors in the normal approximations to these
quantities  depend
upon $\alpha$ and become infinite  as $\alpha \rightarrow 1$. 
Let $h:R\rightarrow R$ 
be a bounded and piecewise continuously differentiable function and
let $\beta = \max(1/2,\alpha)$. 
We use the size-bias version of Stein's method introduced by
Baldi, Rinott and Stein \cite{BRS89} and  prove that, for $n$ large enough,
\[
\Abs{\E h\Bigl(\frac{\inv(\pi)-\mu}{\sigma}\Bigr) - \Phi h} \leq
C\Biggl(\frac{\norm{h}}{\beta(1-\beta)(\beta(1-\beta)n^{1/2}-C_1 n^{-1/2})}
+\frac{\norm{Dh}}{(\beta(1-\beta)n^{1/3}-C_2n^{-2/3})^{3/2}}\Biggr)
\]
where $C$, $C_1$, and $C_2$ are some positive constants, $\Phi h$ is the expectation of $h$ with respect
to the standard normal distribution, and $\mu$ and $\sigma^2$ are the mean
and variance of $\inv(\pi)$, respectively.
If $\alpha \geq 1/2$, then $\beta = \alpha$. Therefore the bound above
diverges as $\alpha \rightarrow 1$.
We prove a similar result for $\des(\pi)$.

Bounds such as the one given in the previous paragraph require $h$ to
be continuous. Goldstein \cite{Goldstein04} has proved a normal approximation
theorem that holds for non-smooth $h$. We use that theorem to prove that
\[\Bigl|\P\Biggl(\frac{\inv(\pi)-\mu}{\sigma} \leq x\Biggr) - \Phi(x)\Bigr|
\leq C(\beta)/\sqrt{n}\]
and that
\[\Bigl|\P\Biggl(\frac{\des(\pi)-\mu}{\sigma} \leq x\Biggr) - \Phi(x)\Bigr|
\leq C(\beta)/\sqrt{n}.\]
These results are contained in Theorems \ref{thm-2-12} and
\ref{thm-2-16}  of this paper. The quantity $C(\beta)$ diverges when
$\alpha \rightarrow 1$. As before $\beta = \max(1/2,\alpha)$.
When the $n$ elements of the multiset are distinct, with $n\geq 2$, we may  use
$\alpha = 1/n$ and $\beta = 1/2$. Therefore the results stated above
imply \eqref{eqn-1-1a} and \eqref{eqn-1-1b}.

The generating function of the number of permutations of a multiset with a given
number of inversions is a rational function. Using this generating function, Diaconis
\cite{DiaconisBook} has shown that the asymptotic distribution of 
$\inv(\pi)$, where $\pi$ is uniformly distributed over permutations of a multiset,
is normal. Theorem \ref{thm-2-12} about $\inv(\pi)$ is accompanied by an error bound 
of the correct order, which is $O(1/\sqrt{n})$, and the dependence of the error bound
on $\alpha$ is also explicitly shown in our theorem. The generating function for
the number of permutations of a multiset with a given number of descents,
which is related to Foata's correspondence, was found
by MacMahon \cite{KnuthVol3} \cite{MacMahonBook}. 
However, normal approximations to this quantity, such as the approximation
given in Theorem \ref{thm-2-16}, do not
seem to be available.

Segments of $\pi(1),\ldots, \pi(n)$ between successive descents, or
runs, are in ascending order. Knuth \cite{KnuthVol3} has stated that
runs are important in the study of sorting algorithms because runs are
segments that are already in sorted order. Among the applications of
descents and inversions to the study of sorting algorithms, multiway
merging with replacement selection merits special mention.  In this
sorting method, the given sequence is first split into runs and the
runs are merged together. Our results are pertinent to sorting
algorithms if the keys used for sorting are allowed to repeat. For
example, Theorem
\ref{thm-2-16} about $\des(\pi)$ gives an idea of how many runs
to expect if multiway merging is used on a sequence of records with repeated
keys.

Descents and inversions have been used as test statistics in the
special case where $\pi$ is a permutation of $\{1,2,\ldots,n\}$. As
already mentioned, early work on runs and descents was stimulated by
statistical applications.  Of the ten empirical tests for the
randomness of a sequence of distinct numbers discussed by Knuth
\cite{KnuthVol2}, one is based on runs and descents.
Taking our results into account, inversions and descents can be used
to test if a given permutation of a multiset of numbers is 
random. There are other ways to test if a given permutation of a
multiset of numbers is random. If a permutation passes $n$ empirical
tests for randomness but fails the $n+1$st, it is not
random. Therefore having a greater number of empirical tests available
makes for more robust testing \cite{KnuthVol2}.

DNA sequences are strings of the four letters $A$, $C$, $G$, and $T$.
It is now well known that these sequences are far from random
\cite{GNC1}.  It has even been suggested that these sequences are
similar to human languages \cite{Mantegna1}.  Some commonly used
compression algorithms such as the Lempel-Ziv method fail to compress
typical DNA sequences however \cite{GNC1}.  The entropy estimates of
DNA sequences given in \cite{GNC1} and \cite{Mantegna1} proceed by
dividing the sequence into blocks in some way. For instance, blocks of
$6$ consecutive letters are considered in \cite{GNC1}. These entropy
estimates show that  DNA sequences are not random.

In Section 3, we report the descents and inversions of the $19$th
chromosome of the human genome mainly as an illustration.  We consider
all $24$ possible orderings of $A$, $C$, $G$, and $T$. With respect to
each of these orderings, a calculation of descents and inversions
shows that the number of descents and inversions of the DNA sequence
departs from the mean by a large multiple of the standard
deviation. It may be of some interest that this method of showing the
DNA sequence to be far from random considers only single letters
without dividing them into blocks.

Although we consider all possible orderings of $A$, $C$, $G$, and $T$,
it must be noted that the molecular weights of the corresponding 
compounds implies the order $C<T<A<G$. This is as natural as any order
one can hope to find among four physical objects.

Our interest in permutations of multisets was provoked by their connection
to riffle shuffles of decks with repeated cards \cite{CV1}.

We do not give explicit numerical constants in our Theorems 
\ref{thm-2-12} and \ref{thm-2-16} about
descents and inversions of permutations of multisets. It is worth
noting that explicit numerical constants are not given for most of the
detailed examples in Stein's book \cite{SteinBook}, and all the
examples in the papers by Baldi, Rinott, and Stein
\cite{BRS89}, by Goldstein and Rinott \cite{GR96}, and by Rinott
and Rotar \cite{RR97}. Furthermore, even the asymptotic normality for
descents of permutations of multisets implied by Theorem
\ref{thm-2-16} is a new result, and so is Theorem \ref{thm-2-12} which
shows the dependence of the bounds for normal approximation on 
the size of the multiset and the parameter $\alpha$ that characterizes
the multiset.

\section{Descents and inversions of permutations of multisets}

If $W\geq 0$ is a non-negative and integrable random variable, the distribution of 
$W^{\ast}$ is said to be $W$-size biased, if $\E(W f(W)) = \E W \E(f(W^\ast))$
for all continuous functions $f$ for which the expectation on the left hand 
side of the equality exists.

Stein's method \cite{Stein72} \cite{SteinBook} refers to the use of auxiliary
randomization to find normal approximations to the distribution of some random
variables. In the theorem below, the
auxiliary randomization requires the construction of $W^{\ast}$ which must be
$W$-size biased. The theorem below can be found in 
\cite{BRS89}, but we follow its formulation in \cite{GR96}.   

\begin{thm}
Let $W$ be a non-negative random variable with $\E W = \mu$ and $\var(W) = \sigma^2$.
Let $W^\ast$ be jointly defined with $W$ such that its distribution is $W$-size
biased. Let $h$ be a function from $R$ to $R$ such that $h$ is continuous
and its derivative $Dh$ is piecewise continuous. Then
\[\Abs{\E h\Bigl(\frac{W-\mu}{\sigma}\Bigr) - \Phi h}
\leq 2 \norm{h} \frac{\mu}{\sigma^2} \sqrt{\var\bigl(\E(W^\ast-W|W)\bigr)}
+\norm{Dh}  \frac{\mu}{\sigma^3}\E(W^\ast-W)^2,\]
where $\Phi h$ is the expectation of $h$ with respect to the
standard normal distribution and $\norm{\cdot}$ is the supremum norm.
\label{thm-2-1}
\end{thm}

When $h$ is the indicator function of the half line $(-\infty,x]$, the following
theorem found in \cite{Goldstein04} applies. Its proof uses a smoothing
inequality and other techniques found in \cite{RR97}.

\begin{thm}
Let $W$ be a non-negative random variable with $\E W = \mu$ and $\var(W) = \sigma^2$.
Let $W^\ast$ be jointly defined with $W$ such that its distribution is $W$-size
biased. Let $\abs{W^\ast-W} \leq B$ and let $A = B/\sigma$. 
Let $B\leq \sigma^{3/2}/\sqrt{6\mu}$.
Then
\[\Abs{P\Bigl(\frac{W-\mu}{\sigma} \leq x\Bigr) - \Phi(x)}
\leq 0.4 A + \frac{\mu}{\sigma}(64A^2+4A^3)+\frac{23\mu}{\sigma^2}
\sqrt{\var\bigl(\E(W^\ast-W|W)\bigr)},\]
where $\Phi$ is the standard normal distribution.
\label{thm-2-2}
\end{thm}

In Theorems \ref{thm-2-1} and \ref{thm-2-2} above, we added the superscript
$*$ to $W$ to denote a random variable with the $W$-size biased
distribution. In the lemma below, random variables $X_i$ with the
superscript $*$ do not necessarily have the $X_i$-size biased
distribution. Here and later, our convention is to use the superscript
$*$ when random variables are constructed as a part of the size biasing
procedure. This notation is due to \cite{BRS89}.

The construction of size biased variables in this paper will be based on the
following lemma found in \cite{BRS89} and \cite{GR96}. 

\begin{lem}
Let $W=X_1+X_2+\ldots+X_n$, where each $X_i$ is a non-negative random variable
with finite mean.
Let $I$ be a random variable which is independent of the $X_i$ and
which satisfies $\P(I=i) = \E X_i/\sum_{j=1}^n\E X_j$. Define $W^\ast$
as $W^\ast = X_1^\ast+ X_2^\ast+\cdots+X_n^\ast$, where for given $I$ $X_I^\ast$ has
the $X_I$-size biased distribution and
\begin{equation}
\P\bigl((X_1^\ast,X_2^\ast,\ldots,X_n^\ast)\in A\bigl|I=i, X_i^\ast=x\bigr)
= \P\bigl((X_1,X_2,\ldots,X_n)\in A\bigl|X_i = x).
\label{eqn-2-1}
\end{equation}
Then $W^\ast$ has the $W$-size biased distribution.
\label{lem-2-3}
\end{lem}
Whenever Lemma \ref{lem-2-3} is applied here, we will find $X_i$ are $0$-$1$ valued
random variables, and the size biased distribution for such variables is 
concentrated at $1$. Therefore, for our purposes, \eqref{eqn-2-1} can be written
as 
\(\P\bigl((X_1^\ast,X_2^\ast,\ldots,X_n^\ast)\in A\bigl|I=i\bigr)
= \P\bigl((X_1,X_2,\ldots,X_n)\in A\bigl|X_i = 1)\).

Let $\pi$ be a uniformly distributed permutation of the multiset
$\{1^{n_1}, 2^{n_2},\ldots,h^{n_h}\}$ and $n=n_1+n_2+\ldots+n_h$. Each $n_a$,
$1\leq a\leq h$, is a positive integer. The symbols $i,j,k,l$, with and without
numerical subscripts, are used to index the set $\{1,2,\ldots,n\}$. The 
symbols $a,b,c,d$ are used to index the set $\{1,2,\ldots,h\}$. We also
assume $n_a \leq \alpha n$ for $1\leq a \leq h$ and for some $\alpha$
in $(0,1)$, $n\geq 4$, and $h\geq 2$.

Define $X_{ij}$, for $i<j$, as $1$ if $\pi(i) > \pi(j)$ and as $0$ otherwise.
Some facts about the joint distribution of $X_{ij}$ will be necessary.
Denote the probabilities 
\begin{align*}
&\P(X_{ij} = 1)\;\text{with $i<j$},\\
&\P(X_{ij_1}=1,X_{ij_2}=1)\; \text{with $i<j_1$ and $i<j_2$},\\
&\P(X_{i_1j}=1,X_{i_2j}=1)\; \text{with $i_1<j$ and $i_2<j$},\\ 
&\P(X_{ik}=1, X_{kj}=1)\; \text{with $i<k<j$,}\\
&\P(X_{i_1j_1}=1,X_{i_2j_2}=1)\;\text{with $i_1<j_1$, $i_2<j_2$, and
$(i_1,j_1)\neq(i_2,j_2)$}
\end{align*}
by $p_1$, $p_2$, $p_3$, $p_4$, and $p_5$, respectively.  Elementary arguments
can be used to deduce formulas, such as $p_1 = \sum_{a<b}n_a n_b/(n(n-1))$
and $p_4= \sum_{a<b<c}n_an_bn_c/(n(n-1)(n-2))$, for $p_1$, $p_2$, $p_3$,
$p_4$, and $p_5$. From such formulas, we deduce
\begin{align}
p_1 &= \frac{n^2-\sum_a n_a^2}{2n(n-1)}\nonumber\\
p_2+p_3+p_4 &= \frac{5n^3/6-n^2+(-3n/2+1)\sum_a n_a^2+(2/3)\sum_a n_a^3}{n(n-1)(n-2)}
\nonumber\\
p_4 &= \frac{n^3/6-(n/2)\sum_a n_a^2 + (1/3)\sum_an_a^3}{n(n-1)(n-2)}\nonumber\\
p_5 &= \frac{n^4/4-n^3+n^2/2+(1/4)\bigl(\sum_a n_a^2\bigr)^2+(-n^2/2+2n-1/2)\sum_{a}n_a^2
      -\sum_a n_a^3}{n(n-1)(n-2)(n-3)}.
\label{eqn-2-2}
\end{align}   
The formulas in \eqref{eqn-2-2} will be used to derive expressions for 
$\var(\inv(\pi))$ and $\var(\des(\pi))$.

The assumption $n_a \leq \alpha n$, for
some $\alpha \in (0,1)$, is used in the two lemmas below. 
The lemmas, however, are worded in terms of $\beta = \max(1/2,\alpha)$
and use the weaker assumption $n_a \leq \beta n$ for
$\beta \in [1/2,1)$.
In both the lemmas the assumption $n_a\leq \beta n$
implies $h\geq 2$.

\begin{lem}
Assume $\beta\in[1/2,1)$, $n_a\geq 0$ for all $a$,
and $\sum_a n_a = n$. If $n_a \leq \beta n$ for $1\leq a\leq h$,
then 
\(2\beta(1-\beta)n^2 \leq n^2-\sum_a n_a^2\leq n^2\) and 
\(3\beta(1-\beta)n^3 \leq n^3-\sum_a n_a^3\leq n^3\).
\label{lem-2-4}
\end{lem}
\begin{proof}
To lower bound $n^2-\sum_a n_a^2$, note that $x\geq y> 0$, $\delta>0$,
and $y-\delta\geq 0$ imply $(x+\delta)^2+(y-\delta)^2 > x^2+y^2$. 
Thus for a given sum $x+y$, the quantity $x^2+y^2$ increases when
the difference $x-y$ is increased.
Thus given $\sum_a n_a = n$ and the constraints
$n_a\geq 0$, the quantity $\sum_a n_a^2$ is
increased whenever two positive numbers are chosen from
$n_a$, $1\leq a \leq h$, and the lesser of them is decreased
and the greater increased by the same amount.
Therefore,
under the constraints $n_a\leq \beta n$,
$\sum_a n_a^2$ is maximum when $n_1=\beta n$,
$n_2=(1-\beta)n$, and $n_a = 0$ for $a>2$. The lower bound for
$n^3-\sum_a n_a^3$ is also obtained when $n_1=\beta n$,
$n_2=(1-\beta)n$, and $n_a = 0$ for $a>2$. The upper bounds are trivial.
\end{proof}

Concerning the lemma below, it is worth noting that $\beta^4-4\beta^4+4\beta-1
= (1-\beta)^2(\beta^2+2\beta-1) > 0$ for $\beta\in [1/2,1)$.
\begin{lem}
Assume $\beta\in[1/2,1)$, $n_a\geq 0$ for all $a$,
and $\sum_a n_a = n$. If $n_a \leq \beta n$ for $1\leq a\leq h$,
\[
(\beta^4-4\beta^2+4\beta-1)n^4 \leq n^4/3+\bigl(\sideset{}{_a}\sum n_a^2\bigr)^2
-(4n/3)\sideset{}{_a}\sum n_a^3 \leq n^4/3.
\]
\label{lem-2-5}
\end{lem}
\begin{proof}
The upper bound follows from the inequality $n\sum_an_a^3 \geq 
\bigl(\sum_a n_a^2\bigr)^2$.

We prove the lower bound assuming $\beta > 1/2$. 
The proof for $\beta=1/2$
can be obtained with minor changes. The proof will make careful use of 
the Kuhn-Tucker conditions as explained in \cite[Theorem 9.2-3]{CiarletBook}.

We attempt to minimize 
$J(n_1,n_2,\ldots, n_h)=\bigl(\sum_an_a^2\bigr)^2-(4n/3)\sum_a n_a^3$ subject to
the affine constraints $\sum_a n_a = n$, $-n_a \leq 0$, and $n_a-\beta n\leq 0$,
where the last two constraints hold for $1\leq a \leq h$.
We assume $n_1\geq n_2\geq\cdots\geq n_h\geq 0$ without loss of generality.

 Let $DJ$ be the gradient vector whose $a$th entry is
\[ \frac{\partial J}{\partial n_a} = 4n_a \sideset{}{_b}\sum n_b^2 - 4n n_a^2 
= 4n_a\sideset{}{_b}\sum n_b(n_b-n_a).\]
The sum of the entries of $DJ$ must be $0$
because the term  $4n_a n_b (n_b-n_a)$ in $\partial J/\partial n_a$
is canceled by the term $4n_a n_b (n_a-n_b)$ in
$\partial J/\partial n_a$.
 If there exists an $a$ such
that $n_1 > n_a > 0$, then the first entry of $DJ$ must be strictly negative\
and therefore some other entry must be strictly positive.

Let $u\in R^h$ be the vector with all entries equal to $1$. Let $v_a\in R^h$
be the vector with its $a$th entry equal to $-1$ and all other entries equal to
$0$. Let $w_1 = -v_1$. Note that $n_a-\beta n=0$ is possible only if $a=1$
as we have assumed $\beta > 1/2$ and  $n_1\geq n_2\geq\cdots$

Suppose $(n_1,n_2,\ldots,n_h)$ is a local minimum of $J$. The
Kuhn-Tucker conditions require that it must be possible to make all
entries of $DJ$ zero by adding multiples of certain vectors. We are
always allowed to add any real multiple of $u$ because the constraint
$\sum_a n_a = n$ is always in force.  We are allowed to add
a positive multiple of $v_a$ if and only if $n_a=0$
because the constraint $-n_a\leq 0$ can then be violated by
making an infinitesimal change to $n_a$.
We are allowed to add a
positive multiple of $w_1$ if and only if $n_1-\beta n = 0$ by
a similar reason.

Let us first consider the type of local minimum where the Kuhn-Tucker conditions
can be satisfied without adding a positive multiple of $w_1$. Suppose
$n_1 > n_a > 0$ for some $a$ for such a local minimum. Then the first entry of
$DJ$ is strictly negative and some other entry is strictly positive. 
If the positive entry is $\partial J/\partial n_b$, then $n_b$ must
be nonzero and therefore $n_b > 0$.
Such
a $DJ$ cannot be made zero by adding a multiple of $u$ and positive multiples
of $v_c$ corresponding to $n_c = 0$.
The only way to make the $b$th entry of $DJ$ equal to $0$ is by
adding a negative multiple of $u$. But this means the first entry
remains negative and nonzero, and the only way to make it $0$ is
by adding a multiple of $w_1$ which is not allowed by assumption.
Therefore any local minimum of this type must have 
$n_1=n_2=\cdots=n_s=n/s$ and $n_c = 0$ for $c > s$, where $2\leq s \leq h$.
The value of $J$ at such a point is $-n^4/(3s^2)$. Since
$s\geq 2$,
\begin{equation}
J \geq -n^4/12
\label{eqn-2-3a}
\end{equation}
at any local minimum of this type. 

We next consider the type of local minimum where it is necessary to add a 
positive multiple of $w_1$ to satisfy the Kuhn-Tucker conditions. At such
a local minimum $n_1 = \beta n$ and $n_1 > n_2\geq \cdots \geq n_h\geq 0$.
Suppose $n_h=0$. Then the first entry of $DJ$ is strictly negative, some other
entry is strictly positive, and the last entry is $0$. We cannot make all those
three entries zero by adding a real multiple of $u$, a positive multiple of $w_1$,
and a positive multiple of $v_h$ to $DJ$. Thus $n_h > 0$.  Next suppose
that the $a$th entry of $DJ$ is not equal to the $b$th entry of $DJ$ for
some $a,b>1$. It is impossible to make the $1$st, $a$th, and $b$th entries of
$DJ$ zero by adding a multiple of $u$ and a positive multiple of $w_1$. Therefore,
all entries of $DJ$ except the first must be equal. The expression for
$\partial J/\partial n_a$ given above is quadratic in $n_a$. Thus we may conclude
that at any local minimum of this type $n_1=\beta n$ and $n_2,\ldots,n_h$
can take on at most two different values. Although the argument assumed 
$h\geq 3$, the conclusion holds when $h=2$ as well. When $h=2$
and a positive multiple of $w_1$ is added to DJ to satisfy
the Kuhn-Tucker conditions, we must have $n_1 = \beta n$
and $n_2 = (1-\beta)n$.

We now consider the value of $J$ assuming that $n_1 = \beta n$, that 
$x$ of the $n_a$s
equal $n_x$, that $y$ of the $n_a$s equal $n_y$, and that $xn_x+yn_y = n(1-\beta)$.
We also assume that $x$ is a positive integer, that $y$ is a non-negative integer, that $x\geq y$, and of course that $n_x$ and $n_y$ are non-negative. Then
\[J = (\beta^4n^4-4\beta^3 n^4/3)+(xn_x^2+yn_y^2)^2 +2\beta^2 n^2(xn_x^2+yn_y^2)
-(4n/3)(xn_x^3+yn_y^3),\]
which we will think of as a sum of four terms. If follows from elementary inequalities
that the minimum of $xn_x^2+yn_y^2$ under the given constraints is 
$n^2(1-\beta)^2/(x+y)$, and that the minimum of $-(4n/3)(xn_x^3+yn_y^3)$
occurs when $n_x = n(1-\beta)$, $x=1$, and $n_y=0$. We can minimize each of 
the four terms of $J$ separately to obtain
\begin{equation}
J \geq \beta^4n^4-4\beta^3 n^4/3-4(1-\beta)^3n^4/3.
\label{eqn-2-3}
\end{equation}
The value of $J$ at any local minimum of the type discussed in the previous paragraph
must either equal or exceed the lower bound in \eqref{eqn-2-3}.

So far, we have proved that the value of $J$ at a local minimum satisfies
the lower bound given by either \eqref{eqn-2-3a} or
\eqref{eqn-2-3}, depending  upon the type of the local minimum.
For $\beta \in [1/2,1)$, $\beta^4-4\beta^3/3-4(1-\beta)^3/3 < -1/12$
by an elementary argument. Therefore the lower bound for $J$ given by
\eqref{eqn-2-3} holds at all local minima and the lower bound for
$J+n^4/3$ stated in the lemma is proved.

\end{proof}

\subsection{Inversions of permutations of multisets}
Let $W = \sum_{i<j} X_{ij}$. Then $W=\inv(\pi)$. We
assume that $\pi$ is uniformly distributed
over permutations of the multiset $\{1^{n_1}, 2^{n_2},\ldots,h^{n_h}\}$.

\begin{lem}
Let $\mu = \E W$ and $\sigma^2 = \var(W)$. Then 
\begin{equation*}
\mu = \frac{n^2 - \sum_a n_a^2}{4}\quad \text{and}\quad
\sigma^2 = \binom{n}{3}\frac{n^5-n^2\sum_a n_a^3}{6n^2(n-1)^2(n-2)}+O(n^2).
\end{equation*}
\label{lem-2-6}
\end{lem}
\begin{proof}
Since $\mu=\binom{n}{2} p_1$, where $p_1 = \E X_{ij}$, and $p_1$ is given
by \eqref{eqn-2-2}, the expression for $\mu$ in the lemma must hold.

We first show that 
\begin{equation}
\sigma^2 = \binom{n}{2}(p_1-p_1^2)+2\binom{n}{3}(p_2+p_3+p_4-3p_1^2)
+6\binom{n}{4}(p_5-p_1^2),
\label{eqn-2-4}
\end{equation}
where the $p_i$ are given by \eqref{eqn-2-2}.
If $\var(W)$ with $W=\sum_{i<j}X_{ij}$ is written as a sum of variances and covariances of the $X_{ij}$,
there are $\binom{n}{2}$ variance terms each of which is equal to 
$p_1-p_1^2$. There are $\binom{n}{3}$ terms of the form
$2\covar(X_{ij_1}, X_{ij_2})$ with $i<j_1<j_2$ and each of those
is equal to $2(p_2-p_1^2)$. We can account for terms of the form
$2\covar(X_{i_1j},X_{i_2j})$ with $i_1<i_2<j$ and of the form
$2\covar(X_{ik}, X_{kj})$ with $i<k<j$ similarly. Thus far we have explained
the first two terms of \eqref{eqn-2-4}. All the other terms in the
expansion of $\var(W)$ are of the form
$2\covar(X_{i_ij_1},X_{i_2j_2})$ with $i_1<j_1$, $i_2<j_2$, and
$(i_1,j_1)<(i_2,j_2)$ in lexicographic order. The last term of  \eqref{eqn-2-4}
follows if we note that the number of such terms is $3\binom{n}{4}$.

The expression for $\sigma^2$ in the lemma is deduced using \eqref{eqn-2-2},
\eqref{eqn-2-4}, and the two inequalities $\sum_a n_a^2 < n^2$ and
$\sum_a n_a^3 < n^3$.
\end{proof}

We now turn to the construction of the size biased variable $W^\ast$ required
by Theorems \ref{thm-2-1} and \ref{thm-2-2}. Let $I$
be uniformly distributed 
over all pairs $(i,j)$ with $1\leq i < j\leq n$
and let it be independent of $\pi$. Let $J=(a,b)$, for
$h\geq a > b\geq 1$, with probability $n_an_b/\sum_{c<d}n_cn_d$, and let
$J$ be independent of both $\pi$ and $I$. Now $\pi^\ast$ is constructed from
$\pi$, $I$, and $J$ as follows. If $I=(i,j)$ and $\pi(i)>\pi(j)$, then
$\pi^\ast=\pi$. If $I=(i,j)$, $\pi(i)\leq \pi(j)$ and $J=(a,b)$, $\pi^\ast$
is constructed in the following steps:
\begin{enumerate}
\item Let $i^\ast$ and $j^\ast$ be uniformly distributed over
the sets $\{i\bigl|\pi(i)=a\}$ and 
$\{j\bigl|\pi(j)=b\}$, respectively. They must be independent
of each other and all other random variables.
\item If $\{i,j\}\cap\{i^\ast,j^\ast\} = \phi$, or $i=i^\ast,j\neq j^\ast$,
or $i\neq i^\ast,j=j^\ast$, exchange $\pi(i)$ with $\pi(i^\ast)$ and
$\pi(j)$ with $\pi(j^\ast)$ to get $\pi^\ast$.
\item If $i=j^\ast,j=i^\ast$, exchange $\pi(i)$ and $\pi(j)$ to get $\pi^\ast$.
\item If $i=j^\ast, j\neq i^\ast$, then 
$\pi^\ast(i) = \pi(i^\ast)$, $\pi^\ast(j) = \pi(j^\ast) = \pi(i)$,
$\pi^\ast(i^\ast) = \pi(j)$, and $\pi^\ast(k) = \pi(k)$ if
$k \neq i, j, i^\ast$.

\item If $i\neq j^\ast, j=i^\ast$, then
$\pi^\ast(i) = \pi(i^\ast) = \pi(j)$, $\pi^\ast(j) = \pi(j^\ast)$, 
$\pi^\ast(j^\ast) = \pi(i)$,  and $\pi^\ast(k) = \pi(k)$ for
$k \neq i, j, j^\ast$.
\end{enumerate}
Finally, $W^\ast = \sum_{i<j} X^\ast_{ij}$, where $X^\ast_{ij}$ is $1$ if
$\pi^\ast(i) > \pi^\ast(j)$ and $0$ otherwise.

We prove below that $W^\ast$ has the $W$-size biased distribution.
If $\pi$ were a  uniformly distributed permutation of 
$\{1,2,\ldots,n\}$, it would be enough to exchange $\pi(i)$
and $\pi(j)$ if $\pi(i)<\pi(j)$ to get $\pi^\ast$. The resulting
$W^\ast$ would have the $W$-size biased distribution. However, since
we are dealing with a multiset here, $\pi(i)=\pi(j)$ is also
a possibility. The construction of $\pi^\ast$ given above is not
as simple mainly because this possibility has to be dealt with.

The following lemma is needed to prove that $W^\ast$ has the $W$-size biased
distribution. Subtraction and union of multisets have the obvious meanings
in the statement of the lemma. The lemma is stated without proof.
\begin{lem}
Let $\pi$ be a uniformly distributed permutation of the multiset
$\{1^{n_1}, 2^{n_2},\ldots,h^{n_h}\}$. If one $a$
out of $n_a$ possible choices is chosen uniformly
from $\pi$ and changed to  $b$, the resulting permutation is a uniformly 
distributed permutation of the multiset $(\{1^{n_1}, 2^{n_2},\ldots,h^{n_h}\}
-\{a\})\cup\{b\}$. 
Similarly, if one of $n_a$ $a$s and one of $n_b$ $b$s are picked uniformly
and independently from $\pi$ and changed to  $c$ and $d$, respectively, then
the resulting permutation is a uniformly distributed permutation of a possibly
new multiset. 
\label{lem-2-7}
\end{lem}

\begin{lem}
The random variable $W^\ast$ has the $W$-size biased distribution.
\label{lem-2-8}
\end{lem}
\begin{proof}
By Lemma \ref{lem-2-3}, it is enough to show that 
\(\P\bigl(\pi^\ast\in A \bigl|I=(i,j)\bigr) = \P\bigl(\pi\in A\bigl| \pi(i)>\pi(j)\bigr)\). Now
\begin{multline*}\P\bigl(\pi^\ast\in A \bigl|I=(i,j)\bigr)
 =
\P\bigl(\pi\in A\bigl|\pi(i)>\pi(j)\bigr)\P(\pi(i)>\pi(j))\\
+\P\bigl(\pi^\ast\in A\bigl|\pi(i)\leq \pi(j), I=(i,j)\bigr)\P(\pi(i)\leq \pi(j)).
\end{multline*}
The first term in the right hand side of the equation above is not
conditioned on $I$ because $\P\bigl(\pi^\ast\in A\bigl| \pi(i)>\pi(j),
I = (i, j)\bigr) = \P\bigl(\pi\in A \bigl| \pi(i) > \pi(j), I = (i,j)\bigr)$,
by the construction of $\pi^\ast$, and because $\pi$ is independent
of $I$.
Thus, if we can show $\P\bigl(\pi^\ast\in A\bigl|\pi(i)\leq \pi(j), I=(i,j)\bigr)
= \P\bigl(\pi\in A\bigl| \pi(i)>\pi(j)\bigr)$, the proof will be complete.

The proof is completed by the sequence of equalities below and the explanation
that follows them.
\begin{align*}
\P&\bigl(\pi^\ast\in A \bigl| \pi(i)\leq \pi(j), I=(i,j)\bigr)\\
=&\sum_{a>b}\P\bigl(\pi^\ast\in A\bigl| \pi(i)\leq\pi(j), I=(i,j), J=(a,b)\bigr)
\P\bigl(\pi(i)=a,\pi(j)=b\bigl| \pi(i)>\pi(j)\bigr)\\
=&\sum_{a>b}\P\bigl(\pi\in A\bigl| \pi(i)=a,\pi(j)=b, I=(i,j), J=(a,b)\bigr)
\P\bigl(\pi(i)=a,\pi(j)=b\bigl| \pi(i)>\pi(j)\bigr)\\
=&\P\bigl(\pi\in A\bigl| \pi(i)>\pi(j)\bigr).
\end{align*}
The first equality is true because $J$ is independent of $\pi$ and $I$,
and $\P\bigl(J=(a,b)\bigr) = P\bigl(\pi(i)=a,\pi(j)=b\bigl|\pi(i)>\pi(j)\bigr)$.
The construction of $\pi^\ast$ from $\pi, I, J$ and Lemma \ref{lem-2-7} imply
the second equality. More specifically, we note that Lemma \ref{lem-2-7}
implies that given $\pi(i)\leq \pi(j)$, $I=(i,j)$ and $J=(a,b)$, the arrangement
$\pi^\ast(1), \pi^\ast(2), \ldots, \pi^\ast(n)$ with the $i$th and the $j$th
numbers {\it struck out} is  a uniformly distributed  permutation of the
multiset $\{1^{n_1},2^{n_2},\ldots,h^{n_h}\}-\{a,b\}$.  
\end{proof}

We now focus on finding a useful upper bound for
 $\var\bigl(\E(W^\ast-W\bigl| \pi)\bigr)$. Given a sequence of
numbers $s_1, s_2,\ldots,s_p$, we  throw $q$ and $r$ into the same set
if and only if $s_q=s_r$. In this way, we get a partition of $\{1,2,\ldots,p\}$ into
sets, and we may arrange the sets of the partition so that the values of $s_q$
for $q$ in the set increase.
We refer to such an ordered partition of $\{1,2,\ldots,p\}$ as the {\it relative
order} of $s_1,s_2,\ldots,s_p$. For our purpose, it is sufficient to note that
the number of possible relative orders is bounded by $2^p p!$.

\begin{lem}
Let $P_1$ be the probability that $\pi(1),\pi(2), \ldots,\pi(p)$ occur
in a certain relative order when $\pi$ is a uniformly distributed permutation
of $\{1^{n_1}, 2^{n_2},\ldots, h^{n_h}\}$, and let that probability be $P_2$
if $\pi$ is a uniformly distributed permutation of the multiset
$\{1^{n'_1},2^{n'_2},\ldots, h^{n'_h}\}$. Assume that $n_a\geq n'_a$,
$\sum_a (n_a-n'_a)\leq 5$. We allow $n'_a=0$. If $p\leq 5$ then
\(\Abs{P_1-P_2} \leq C/n\)
for some constant $C$.
\label{lem-2-9}
\end{lem}
\begin{proof}
The proof is obtained by writing down formulas for $P_1$ and $P_2$. We show 
the proof for the relative order $\pi(1)<\pi(2)<\cdots<\pi(p)$.

Let $n'=\sum_an'_a$. The probability $P_1$ is given by
\begin{equation}
\frac{\sum n_{a_1}n_{a_2}\ldots n_{a_p}}{n(n-1)\ldots(n-p+1)},
\label{eqn-2-5}\
\end{equation}
where the sum is taken over $1\leq a_1<a_2<\cdots<a_p\leq h$. The formula
for $P_2$ is obtained by adding a prime to all the $n$s in \eqref{eqn-2-5}.
Now
\[ P_1 - P_2 = 
\frac{\sum n_{a_1}n_{a_2}\ldots n_{a_p}- n'_{a_1}n'_{a_2}\ldots n'_{a_p}}
{n(n-1)\ldots(n-p+1)} - 
P_2 \Biggl(1 -\frac{n'(n'-1)\ldots(n'-p+1)}{n(n-1)\ldots(n-p+1)}\Biggr),\]
$0\leq n-n'\leq 5$, and
\[n_{a_1}n_{a_2}\ldots n_{a_p}- n'_{a_1}n'_{a_2}\ldots n'_{a_p}
\leq n_{a_1}n_{a_2}\ldots n_{a_p}((n_{a_1}-n'_{a_1})/n_{a_1}+\cdots+
(n_{a_p}-n'_{a_p})/n_{a_p}),\]
together imply
$\abs{P_1-P_2} \leq C/n$.
\end{proof}

\begin{lem}
Let $f\bigl(\pi(1),\pi(2),\ldots,\pi(p)\bigr)$ and $g\bigl(\pi(p+1),\pi(p+2),\ldots,\pi(p+q)\bigr)$ be functions that depend
only upon the relative order of their argument lists. Assume that $\abs{f}$, $\abs{g}$,
$p$, and $q$ are all upper bounded by  $5$. If $\pi$ is
a uniformly distributed permutation of the multiset
$\{1^{n_1},2^{n_2},\ldots,h^{n_h}\}$, then
\[\Bigl|\covar\bigl(f\bigl(\pi(1),\pi(2),\ldots,\pi(p)\bigr),
g\bigl(\pi(p+1),\pi(p+2),\ldots,\pi(p+q)\bigr)\bigr)\Bigr| \leq C/n\]
for some constant $C$.
\label{lem-2-10}
\end{lem}
\begin{proof}
It is enough to consider $f$ and $g$ to be indicator functions that are $1$
for a certain relative order of their argument lists and $0$ for all other 
relative orders. All other $f$ and $g$ are linear combinations of a constant 
number of indicator functions with coefficients that are bounded by constants.

We state the proof assuming $f$ and $g$ are $1$ if their arguments are in 
strictly increasing order and $0$ otherwise. Let $\P(f=1)=P_1$ and
$\P(g=1)=P_2$. Then
\begin{multline*}
\P(fg = 1) = \sum
\P\bigl(\pi(1)<\pi(2)<\cdots<\pi(p)\bigl|\pi(p+1)=a_1,\ldots,\pi(p+q)=a_q\bigr)\\
\P\bigl(\pi(p+1)=a_1,\ldots,\pi(p+q)=a_q\bigr),
\end{multline*}
where the sum is over $1\leq a_1<a_2<\cdots<a_q\leq h$.
By the previous Lemma \ref{lem-2-9}, each conditional probability in the sum above
is $P_1 + O(1/n)$. Therefore, $\P(fg=1)=P_1P_2 + O(1/n)$ and
$\covar(f,g)=O(1/n)$.
\end{proof}

\begin{lem}
\[\var\bigl(\E\bigl(W^\ast-W\bigl|\pi\bigr)\bigr) \leq \frac{Cn^5}{(n^2-\sum_an_a^2)^2}\]
for some constant $C$.
\label{lem-2-11}
\end{lem}
\begin{proof}
If $\pi(i) > \pi(j)$, $\E\bigl(W^\ast-W\bigl|\pi, I=(i,j)\bigr) = 0$.
If $\pi(i)\leq \pi(j)$,
\[ \E\bigl(W^\ast-W\bigl|\pi, I=(i,j)\bigr) =
\frac{1}{\sum_{a>b}n_an_b} \sum_{i^\ast,j^\ast}\sum_{l=1}^n \psi_\pi(i,j,i^\ast,j^\ast,l),\]
where $(i^\ast,j^\ast)$ takes all $\sum_{a>b}n_an_b$ possible values with
$\pi(i^\ast) > \pi(j^\ast)$ and $\psi_\pi(i,j,i^\ast,j^\ast,l)$ is the change
in the number of inversions between position $l$ and positions $i,j,i^\ast,j^\ast$
when $\pi(i)$, $\pi(j)$, $\pi(i^\ast)$, $\pi(j^\ast)$ are exchanged to construct
$\pi^\ast$. Note that $\abs{\psi_\pi} \leq 4$. We now have
\begin{equation}
\E\bigl(W^\ast-W\bigl|\pi\bigr) 
= \frac{1}{\binom{n}{2}\sum_{a>b}n_an_b}\sum_{i,j}\sum_{i^\ast,j^\ast}\sum_{l=1}^n
\psi_\pi(i,j,i^\ast,j^\ast, l),
\label{eqn-2-6}
\end{equation}
where $i,j$ take all values satisfying $1\leq i < j\leq h$ and $\pi(i)\leq \pi(j)$,
and where $i^\ast,j^\ast$ take values as already indicated. 

We use \eqref{eqn-2-6} to write
$\var\bigl(\E\bigl(W^\ast-W\bigl|\pi\bigr)\bigr)$ as a sum of variance and
covariance terms. The number of variance terms is bounded by $n^5$. The number
of covariance terms
\begin{equation}
\covar\bigl(\psi_\pi(i_1,j_1,i_1^\ast,j_1^\ast,l_1), \psi_\pi(i_2,j_2,i_2^\ast,j_2^\ast,l_2)\bigr)
\label{eqn-2-7}
\end{equation}
with $\{i_1,j_1,i_1^\ast,j_1^\ast,l_1\}\cap\{i_2,j_2,i_2^\ast,j_2^\ast,l_2\}
\neq \phi$ is fewer than $25n^9$. Since $\abs{\psi_\pi} \leq 4$, the contribution
of the variance terms and covariance terms with the property just described
is bounded by $16(n^5+25 n^9)/\bigl(\binom{n}{2}\frac{1}{2}(n^2-\sum_an_a^2)\bigr)^2$. We have used $\sum_{a>b} n_a n_b = \frac{1}{2}(n^2-\sum_a n_a^2)$
to obtain this bound.

Covariance terms of the form \eqref{eqn-2-7} with 
\(\{i_1,j_1,i_1^\ast,j_1^\ast,l_1\}\cap\{i_2,j_2,i_2^\ast,j_2^\ast,l_2\}
= \phi\)
remain to be considered. The number of such terms is fewer than $n^{10}$. 
Lemma \ref{lem-2-10} can be applied to argue that such covariances are 
$O(1/n)$ as we may use the fact that $\pi$ is uniformly distributed to
assume $i_1,j_1,i_1^\ast,j_1^\ast,l_1=1,2,3,4,5$ and
$i_2,j_2,i_2^\ast,j_2^\ast,l_2=6,7,8,9,10$ with no loss of generality.
The proof can now be easily completed.
\end{proof}

\begin{thm} Let $\pi$ be a uniformly distributed permutation of the multiset
$\{1^{n_1}, 2^{n_2},\ldots, h^{n_h}\}$, where $n_a\in Z^+$ for $1\leq a \leq h$.
Assume that $\alpha\in (0,1)$ is fixed and that $n_a\leq \alpha n$ for
$1\leq a \leq h$. 
Let $\beta = \max(1/2,\alpha)$.
Let $h:R\rightarrow R$ be a bounded continuous function with 
bounded piecewise continuous derivative $Dh$. Then for $n > n_0(\beta)$,
\[
\Abs{\E h\Bigl(\frac{\inv(\pi)-\mu}{\sigma}\Bigr) - \Phi h} \leq
C\Biggl(\frac{\norm{h}}{\beta(1-\beta)(\beta(1-\beta)n^{1/2}-C_1 n^{-1/2})}
+\frac{\norm{Dh}}{(\beta(1-\beta)n^{1/3}-C_2n^{-2/3})^{3/2}}\Biggr)
\]
where $C$, $C_1$, and $C_2$ are some positive constants, $\Phi h$ is the expectation of $h$ with respect
to the standard normal distribution, and $\mu$ and $\sigma^2$ are the mean
and variance of $\inv(\pi)$, respectively.

If $C(\beta)$ is allowed to depend upon $\beta$, we may assert
\[\Bigl|\P\Biggl(\frac{\inv(\pi)-\mu}{\sigma} \leq x\Biggr) - \Phi(x)\Bigr|
\leq C(\beta)/\sqrt{n}\]
for some positive constant $C(\beta)$.
\label{thm-2-12}
\end{thm}
\begin{proof}
Let $W=\inv(\pi)$. By Lemmas \ref{lem-2-4} and \ref{lem-2-6},
$\sigma^2 \geq (\beta(1-\beta)/12)n^3+O(n^2)$ and
$\mu\leq n^2/4$. By Lemmas \ref{lem-2-4} and \ref{lem-2-11},
$\var\bigl(\E\bigl(W^\ast-W\bigl|\pi\bigr)\bigr)\leq Cn/(\beta(1-\beta))^2$
for some constant $C$.
By construction of the size biased variable $W^\ast$, $\Abs{W^\ast-W}
\leq 4n$, and therefore $\E\bigl(W^\ast-W\bigr)^2\leq 16n^2$. If we note
that $\var\bigl(\E\bigl(W^\ast-W\bigl|W\bigr)\bigr)\leq \var\bigl(\E\bigl(W^\ast-W\bigl|\pi\bigr)\bigr)$, Theorem \ref{thm-2-1}
can be applied to prove the first part of this theorem.

 The second part 
is proved using Theorem \ref{thm-2-2}. By construction of $W^\ast$,
$\abs{W^\ast - W} \leq 4n$. Therefore we can take $B=4n$. The 
inequality $B \leq \sigma^{3/2}/\sqrt{6\mu}$ must  hold
for large enough $n$  by bounds
for $\sigma $ and $\mu$ given above.
\end{proof}

\subsection{Descents of permutations of multisets}
Let $W=X_{12}+X_{23}+\cdots+X_{n-1,n}$. Then $W=\des(\pi)$, with $\pi$ uniformly
distributed over permutations of the multiset $\{1^{n_1},2^{n_2},\ldots,h^{n_h}\}$.

\begin{lem}
Let $\mu = \E W$ and $\sigma^2=\var(W)$. Then
\begin{equation*}
\mu = \frac{n^2-\sum_a n_a^2}{2n}\quad \text{and}\quad
\sigma^2 = \frac{n^4/3+\Bigl(\sum_a n_a^2\Bigr)^2-(4n/3)\sum_a n_a^3}{4n(n-1)^2}
+ O(1).
\end{equation*}
\label{lem-2-13}
\end{lem}
\begin{proof}
Since $\mu=(n-1)p_1$, where $p_1 = \E X_{ij}$, and $p_1$ is given
by \eqref{eqn-2-2}, the expression for $\mu$ in the lemma must hold.

We first show that 
\begin{equation}
\sigma^2 = (n-1)(p_1-p_1^2)+2(n-2)(p_4-p_1^2)+(n-2)(n-3)(p_5-p_1^2),
\label{eqn-2-8}
\end{equation}
where the $p_i$ are given by \eqref{eqn-2-2}. If $\var(W)$, with 
$W=X_{12}+X_{23}+\cdots+X_{n-1,n}$, is written as the sum
of variances and covariances of the $X_{i,i+1}$,
there are $(n-1)$ variance terms, each equal to $p_1-p_1^2$. There
are $(n-2)$ covariance terms of the form $\covar\bigl(X_{i,i+1},X_{i+1,i+2}\bigr)$
each equal to $p_4-p_1^2$. The remaining covariance terms are all equal to
$p_5-p_1^2$.

The expression for $\sigma^2$ in the lemma is deduced using \eqref{eqn-2-2},
\eqref{eqn-2-8}, and the two inequalities $\sum_a n_a^2 < n^2$ and
$\sum_a n_a^3 < n^3$.
\end{proof}

The construction of the size biased variable $W^\ast$ is the same as the construction
for inversions given immediately after Lemma \ref{lem-2-6} with the following 
differences. The random variable $I$ must be equal to one of $(1,2),(2,3),\ldots,
(n-1,n)$ with equal probability. In the construction of $\pi^\ast$, the symbol
$j$ must be replaced everywhere by $i+1$. Finally, 
$W^\ast = X^\ast_{12}+X^\ast_{23}+\cdots+X^\ast_{n,n-1}$, where 
$X^\ast_{ij}$ is $1$ if $\pi^\ast(i)>\pi^\ast(j)$ and $0$ otherwise.

\begin{lem}
The random variable $W^\ast$ has the $W$-size biased distribution.
\label{lem-2-14}
\end{lem}
\begin{proof}
Similar to the proof of Lemma \ref{lem-2-8}.
\end{proof}

\begin{lem}
\[\var\bigl(\E\bigl(W^\ast-W\bigl|\pi\bigr)\bigr)\leq 
\frac{Cn^5}{n^2(n^2-\sum_a n_a^2)^2}\]
for some constant $C$.
\label{lem-2-15}
\end{lem}
\begin{proof}
By arguing as in the proof of Lemma \ref{lem-2-11}, we get
\begin{equation}
\E\bigl(W^\ast-W\bigl|\pi\bigr) = \frac{2}{n(n^2-\sum_a n_a^2)} \sum_i
\sum_{i^\ast, j^\ast} \psi_\pi(i,i^\ast, j^\ast).
\label{eqn-2-9}
\end{equation}
In \eqref{eqn-2-9}, $i$ takes all values such that $\pi(i)\leq \pi(i+1)$, $(i^\ast,j^\ast)$
takes all values such that $\pi(i^\ast) > \pi(j^\ast)$, and
$\psi_\pi(i,i^\ast, j^\ast) = \des(\pi^\ast)-\des(\pi)$, where $\pi^\ast$ is
constructed by exchanging $\pi(i), \pi(i+1), \pi(i^\ast), \pi(j^\ast)$ as
described. Note that $\abs{\psi_\pi} \leq 7$.

We use \eqref{eqn-2-9}
to write $\var\bigl(\E\bigl(W^\ast-W\bigl|\pi\bigr)\bigr)$ as the sum of
variance and covariance terms. There are $O(n^3)$ variance terms of the
$\var(\psi_\pi)$. The number of terms of the form
\begin{equation}
\covar\bigl(\psi_\pi(i_1,i^\ast_1,j^\ast_1),\psi_\pi(i_2,i^\ast_2,j^\ast_2)\bigr), 
\label{eqn-2-10}
\end{equation}
where one of the numbers $\{i_1,i^\ast_1,j^\ast_1\}$ differs from one of the 
numbers $\{i_2,i^\ast_2,j^\ast_2\}$ by $3$ or less in magnitude is 
$O(n^5)$. The magnitude of such covariance terms and of the variance terms
is bounded by $49$.  The number of covariance terms of the form \eqref{eqn-2-10}
where none of the numbers $\{i_1,i^\ast_1,j^\ast_1\}$ differs from any one of the 
numbers $\{i_2,i^\ast_2,j^\ast_2\}$ by $3$ or less in magnitude is 
$O(n^6)$. By Lemma \ref{lem-2-10}, the magnitude of such covariance terms is
$O(1/n)$. The proof is now easily completed.
\end{proof}

It is worth noting again that
$\beta^4-4\beta^2+4\beta-1 = (1-\beta)^2(\beta^2+2\beta-1)>0$
for $\beta \in [1/2,1)$.

\begin{thm} Let $\pi$ be a uniformly distributed permutation of the multiset
$\{1^{n_1}, 2^{n_2},\ldots, h^{n_h}\}$, where $n_a\in Z^+$ for $1\leq a \leq h$.
Assume that $\alpha\in(0,1)$ is fixed and that $n_a\leq \alpha n$ for
$1\leq a \leq h$.
Let $\beta = \max(1/2,\alpha)$.
Let $h:R\rightarrow R$ be a bounded continuous function with 
bounded piecewise continuous derivative $Dh$. Then for $n > n_0(\beta)$,
\begin{multline*}
\Abs{\E h\Bigl(\frac{\des(\pi)-\mu}{\sigma}\Bigr) - \Phi h} \leq
C\Biggl(\frac{\norm{h}}
{\beta(1-\beta)\sqrt{n}(\beta^4-4\beta^2+4\beta-1-C_1n^{-1})}\\
+\frac{\norm{Dh}}{((\beta^4-4\beta^2+4\beta-1)n^{1/3}-C_2n^{-2/3})^{3/2}}\Biggr)
\end{multline*}
where $C$, $C_1$, and $C_2$ are some positive constants, $\Phi h$ is the expectation of $h$ with respect
to the standard normal distribution, and $\mu$ and $\sigma^2$ are the mean
and variance of $\des(\pi)$, respectively.

If $C(\beta)$ is allowed to depend upon $\beta$, we may assert
\[\Bigl|\P\Biggl(\frac{\des(\pi)-\mu}{\sigma} \leq x\Biggr) - \Phi(x)\Bigr|
\leq C(\beta)/\sqrt{n}\]
for some positive constant $C(\beta)$.
\label{thm-2-16}
\end{thm}
\begin{proof}
Let $W=\des(\pi)$. By Lemmas \ref{lem-2-4}, \ref{lem-2-5} and
\ref{lem-2-13}
$\sigma^2 \geq ((\beta^4-4\beta^2+4\beta-1)/4)n+O(1)$ and
$\mu\leq n/2$. By Lemmas \ref{lem-2-4} and \ref{lem-2-15},
$\var\bigl(\E\bigl(W^\ast-W\bigl|\pi\bigr)\bigr)\leq 
C/\bigl(n\beta^2(1-\beta)^2\bigr)$ for some constant $C$.
By construction of the size biased variable $W^\ast$, $\Abs{W^\ast-W}
\leq 7$, and therefore $\E\bigl(W^\ast-W\bigr)^2\leq 49$. If we note
that $\var\bigl(\E\bigl(W^\ast-W\bigl|W\bigr)\bigr)\leq \var\bigl(\E\bigl(W^\ast-W\bigl|\pi\bigr)\bigr)$, Theorem \ref{thm-2-1}
can be applied to prove the first part of this theorem.

 The second part 
is proved using Theorem \ref{thm-2-2}. By construction of $W^\ast$,
$\abs{W^\ast - W} \leq 8$. Therefore we can take $B=8$. The 
inequality $B \leq \sigma^{3/2}/\sqrt{6\mu}$ must  hold
for large enough $n$  by bounds
for $\sigma $ and $\mu$ given above.

\end{proof}

\section{Descents and inversions of the human genome}

The human genome consists of 24 chromosomes, each of which is a sequence of bases
labeled A, C, G, or T. The $19$th chromosome has the following counts for 
the four bases (see \cite{IHGC}):
\[n_A = 14383026\quad n_C=13473774\quad n_G=13506612\quad n_T=14422243.\]
The version of the human genome reported in \cite{IHGC} has $341$ gaps. The
$19$th chromosome has only three gaps in the middle. We ignored these gaps when counting the number of inversions and descents.  

\begin{table}
\begin{center}
\begin{tabular}{|c||r|r|r|r|}\hline
{}& A & C & G & T \\ \hline\hline 
A  & 4229414 & 2833985 & 4154304 & 3165323\\ \hline  
C  & 4221129 & 4044958 & 1057112 & 4150574\\ \hline 
G  & 3423863 & 3180474 & 4056078 & 2846197\\ \hline  
T  & 2508620 & 3414357 & 4239118 & 4260148\\ \hline  
\end{tabular}
\end{center}

\begin{center}
\begin{tabular}{|c||r|r|r|r|}\hline
{} & A & C & G & T \\\hline\hline 
A & 103435711266825  & 94175991781325  & 94404662110136  & 103982892949612\\ \hline  
C &  99617649978799  & 90771286164651  & 90984870248490  & 100143945584446\\ \hline  
G &  99861289457776  & 91000167345198  & 91214277105966  & 100388853252680\\ \hline  
T & 103452603097706  & 94178097170636  & 94406787118036  & 104000539364403\\ \hline 
\end{tabular}
\end{center}
\caption[xyz]{The first table above reports the number of occurrences
of $\pi(i)=x$ and $\pi(i+1)=y$.
The second table reports the number of occurrences of $\pi(i)=x$
and $\pi(j)=y$, with $i<j$.
The permutation $\pi$ corresponds to chromosome $19$, and $x$ and $y$
can be A, C, G, or T.}
\label{tbl-3-1}
\end{table}

\begin{table}
\begin{center}
\begin{tabular}{||c|r|r||c|r|r||}\hline
Order & $(\des-\mu_d)/\sigma_d$& $(\inv-\mu_i)/\sigma_i$ & 
Order & $(\des-\mu_d)/\sigma_d$& $(\inv-\mu_i)/\sigma_i$ \\\hline\hline
A,C,G,T  &   $36.13$  &      $-11.64$ & 
C,A,G,T  & $-628.47$   &     $-92.58$ \\\hline
G,A,C,T  & $-631.23$   &     $-93.03$&
A,G,C,T  & $-981.20$  &      $-11.86$\\\hline
C,G,A,T  & $-278.50$ & $-173.74$ &  
G,C,A,T  & $-1295.83$  &     $-173.96$\\ \hline
T,C,A,G  &  $-628.47$  &       $93.03$&
C,T,A,G  &  $-981.20$  &        $4.29$\\\hline
A,T,C,G  &  $-278.50$  &      $166.08$&
T,A,C,G  &  $  36.13$  &      $173.96$\\\hline
C,A,T,G  & $-1295.83$  &       $-3.60$&
A,C,T,G  & $-631.23$  &       $77.34$\\\hline
A,G,T,C  & $-628.47$  &       $76.87$&
G,A,T,C  &  $-278.50$  &       $-4.29$\\\hline
T,A,G,C  &  $-981.20$  &      $173.74$&
A,T,G,C  & $-1295.83$  &      $165.85$\\\hline
G,T,A,C  &    $36.13$  &        $3.60$&
T,G,A,C  &  $-631.23$  &       $92.58$\\ \hline
T,G,C,A  & $-1295.83$  &       $11.64$&
G,T,C,A  &  $-628.47$  &      $-77.34$\\\hline
C,T,G,A  & $-631.23$   &     $-76.87$&
T,C,G,A  &  $-278.50$  &       $11.86$\\\hline
G,C,T,A  &  $-981.20$  &     $-166.08$&
C,G,T,A  &    $36.13$  &     $-165.85$\\ \hline
\end{tabular}
\end{center}
\caption[xyz]{This table 
reports the normalized
number of descents and inversions of the 
$19$th chromosome, when the orders shown in the first and the fourth columns
are considered increasing.} 
\label{tbl-3-2}
\end{table}

From Lemmas \ref{lem-2-6} and \ref{lem-2-13}, and their proofs, we find
the expected number of descents and inversions to be
$\mu_d = 2.0912146861\times 10^{7}$ and $\mu_i = 5.8329890505\times 10^{14}$,
respectively. The standard deviations are
$\sigma_d = 2.0871959423\times 10^3$ and $\sigma_i = 6.7231321079\times 10^{10}$.
Data about the $19$th chromosome reported in Table \ref{tbl-3-1} 
can be used to
calculate the number of descents and inversions for any ordering of 
A, C, G, and T. By Theorems \ref{thm-2-12} and \ref{thm-2-16}, the number of descents
and inversions must have a distribution that is close to the normal 
distribution if $\pi$ is
a uniformly distributed permutation of the bases in the $19$th chromosome. The
number of descents and inversions in the $19$th chromosome itself is reported
in Table \ref{tbl-3-2}
for all possible orderings of  $A$, $C$, $G$, and $T$
and with suitable normalization. From each line of this table, we may infer that the null hypothesis
stating the $19$th chromosome to be a random permutation of its bases is very
unlikely to hold.

Estimations of the entropy of DNA sequences can be found in
\cite{GNC1} and \cite{Mantegna1}. Those estimates too imply that DNA
sequences are far from random. We note that Table \ref{tbl-3-2}
assumes the number of $A$s, $C$s, $G$s, and $T$s to be given and
computes a statistic to test if their arrangement in a sequence is
random.  This is a different notion of randomness from that of
entropy. For instance, it is possible for a sequence to have $A$ for
$90\%$ of its letters which would mean that the sequence can be
significantly compressed. Yet the arrangement of the letters could be
generated randomly.

In the bounds given by Theorems \ref{thm-2-12} and \ref{thm-2-16} the
constants $C(\beta)$ are not determined explicitly. In this example
$n$ is greater than $5 \times 10^7$. For the large departures from the
mean that are seen in Table \ref{tbl-3-2}, it is reasonable to assume
that the probabilities of finding such departures, if the sequence
were a uniformly distributed permutation, are less than $.001$.  Such
a bound is implied in most cases by Chebyshev's inequality.  Yet even
this is surely an overestimate.  For uniformly distributed $\pi$, the
probabilities that $\des(\pi)$ and $\inv(\pi)$ depart from their means
by a certain amount appear to fall off at least as fast as the bell
curve does away from zero. Therefore, for large deviations from the
mean, the bounds given by Theorems
\ref{thm-2-12} and \ref{thm-2-16} are not accurate and better bounds
would be desirable.

\section{Acknowledgments}

The authors thank D. Burns and J. Fulman for helpful discussions, and
the referees for their work.


\bibliography{references}
\bibliographystyle{plain}
\end{document}